\documentclass[aapm,graphicx,article,superscriptaddress]{revtex4-1}
\usepackage{subfigure}
\usepackage[pdftex]{graphicx}
\usepackage[T1]{fontenc}
\usepackage[utf8]{inputenc}
\usepackage[colorlinks=true,linkcolor=blue,citecolor=blue]{hyperref}
\DeclareGraphicsExtensions{.pdf,.jpeg,.png}
\usepackage{tikz}
\usetikzlibrary{arrows}
\usepackage{float}
\usepackage{enumitem}
\usepackage{mathrsfs}
\usepackage{mathtools}
\usepackage{amsfonts}
\usepackage{amssymb}
\usepackage{amsmath}
\usepackage[mathlines]{lineno}
\usepackage{titlesec}
\usepackage{multirow}
\usepackage{rotating}

\begin{document}


\title{Towards Robust Adaptive Radiation Therapy Strategies} 



\author{Michelle B\"ock}
\affiliation{KTH Royal Institute of Technology, Stockholm, Sweden}
\affiliation{RaySearch Laboratories AB, Stockholm, Sweden}
\author{Kjell Eriksson}
\affiliation{RaySearch Laboratories AB, Stockholm, Sweden}
\author{Anders Forsgren}
\affiliation{KTH Royal Institute of Technology, Stockholm, Sweden}
\author{Bj\"orn H\aa rdemark}
\affiliation{RaySearch Laboratories AB, Stockholm, Sweden}


\date{\today}

\begin{abstract}
\begin{description}
\item[Purpose] To set up a framework combining robust treatment planning with adaptive reoptimization in order to maintain high treatment quality, to respond to interfractional geometric variations and to identify those patients who will benefit the most from an adaptive fractionation schedule.
\item[Methods] The authors propose robust adaptive strategies based on stochastic minimax optimization for a series of simulated treatments on a one-dimensional patient phantom. The plan applied during the first fractions should be able to handle anticipated systematic and random errors. Information on the individual geometric variations is gathered at each fraction. At scheduled fractions, the impact of the measured errors on the delivered dose distribution is evaluated. For a patient having received a dose that does not satisfy specified plan quality criteria, the plan is reoptimized based on these individually measured errors. The reoptimized plan is then applied during subsequent fractions until a new scheduled adaptation becomes necessary. In this study, three different adaptive  strategies are introduced and investigated. (i) In the first adaptive strategy, the measured systematic and random error scenarios and their assigned probabilities are updated to guide the robust reoptimization. (ii) In the second strategy, the grade of conservativeness is adapted in response to the measured dose delivery errors. (iii) In the third strategy, the uncertainty margins around the target are recalculated based on the measured errors. The simulated treatments are subjected to systematic and random errors that are either similar to the anticipated errors or unpredictably larger in order to critically evaluate the performance of these three adaptive strategies.
\item[Results] According to the simulations, robustly optimized treatment plans provide sufficient treatment quality for those treatment error scenarios similar to the anticipated error scenarios. Moreover, combining robust planning with adaptation leads to improved organ at risk protection. In case of unpredictably larger treatment errors, the first strategy in combination with at most weekly adaptation performs best at notably improving treatment quality in terms of target coverage and OAR protection in comparison with a non-adaptive approach and the other adaptive strategies.
\item[Conclusion] The authors present a framework that provides robust plan reoptimization or margin adaptation of a treatment plan in response to interfractional geometric errors throughout the fractionated treatment. According to the simulations, these robust adaptive treatment strategies are able to identify candidates for an adaptive treatment, thus giving the opportunity to provide individualized plans, and improve their treatment quality through adaptation. The simulated robust adaptive framework is a guide for further development of optimally controlled robust adaptive therapy models.

\end{description}

\end{abstract}

\maketitle 

\section{Introduction}
\label{sec:intro}
Adaptive radiation therapy (ART) is an approach that allows for adaptation during the course of treatment to patient-specific anatomical and/or biological changes that could not have been accounted for during the planning process.  By combining state-of-the-art medical imaging and treatment plan optimization to precisely deliver the optimized dose plan to the patient, highly individualized radiation therapy for cancer patients can become a reality. Treatment plan optimization software is supposed to respond to geometric variations in order to mitigate deviations of the delivered dose from the intended dose distribution. 

The robust approach to radiation therapy planning includes uncertainties in the treatment delivery process. Uncertainties are caused for example by internal organ motion, geometric changes in the patient's anatomy and setup errors, which can be included in the optimization problem in different ways. The framework presented in this work considers interfractional geometric variations only. Probabilistic methods optimize the expectation value of the objective function by incorporating all possible error scenarios. The probability of occurrence for each error scenario is modeled as a weighting factor in the objective function \cite{unkelbach2004inclusion}. Mini-max optimization is another method used to achieve a treatment plan satisfying specified plan quality criteria despite all errors that may occur during fractionation. The goal here is to find the optimal treatment plan that will guarantee a successful treatment in presence of the worst case scenario \cite{Albin}. Inclusion of the worst case scenario may result in a robust plan in terms of target coverage, but may also result in a too high dose to the adjacent critical structures. Therefore it is desirable to be also robust in terms of healthy tissue sparing.
The motivation behind combining robust optimization with adaptive strategies is to reduce the negative effects of uncertainties on healthy tissue close to the target during the fractionated treatment. By using a feedback process during the fractionation schedule the uncertainties in treatment delivery can be measured and their impact on the cumulative dose distribution can be quantified. This may allow for adaptation of the robust optimization problem. Three different adaptive strategies are applied to and evaluated on a one-dimensional phantom. The evaluations of the delivered dose are carried out according to predefined schedules. The treatment plan, if necessary, is reoptimized for the subsequent fractions in accordance with the three adaptive strategies. Moreover, measures are introduced to identify those patients from the simulations who will benefit the most from an adaptive treatment approach, since previous studies have concluded that not all patients may benefit to the same extent from an adaptive treatment~\cite{olteanu2014comparative,schwartz2013adaptive}. In the first strategy, the measured treatment errors are evaluated in order to update the uncertainty scenarios of the systematic and random errors and their probabilities in the robust optimization problem. The adapted plan is reoptimized with respect to the individual uncertainty distribution. In the second strategy, the grade of conservativeness is adapted depending on the extent of the geometric errors and their impact on the delivered dose distribution. In the third strategy, the safety margin around the clinical target volume is recalculated based on the measured systematic and random errors giving an adapted planning target volume in the optimization problem. By studying an idealized geometry it is possible to focus on the mathematical concepts of robust adaptive strategies.

The aim of this work is to introduce a proof-of-concept framework of an adaptive treatment planning approach. By studying the effect of the three strategies combined with various adaptation schedules on the cumulative dose distribution and treatment quality, in the presence of geometric uncertainties, the group benefiting the most from a robust adaptive approach may be identified and which schedule and strategy combination may be most suitable for this particular group of patients. 
\section{Geometry Model}\label{sec:Geometry}For this study a one-dimensional phantom as depicted below in Figure \ref{fig:Geometry} is considered. The model schematically represents a slice of a two-dimensional phantom or the intersection of a sagittal and transversal cut of a three-dimensional geometry. We consider to irradiate the target with a perpendicular oriented field. The phantom is descretized in 101 voxels and the absorbed dose in each voxel is modeled by Gaussian functions at a spacing of 1 mm with a standard deviation of 3 mm. Acquiring the planning CT (pCT) is the first step in the radiation treatment planning process. Based on this image a radiation oncologist defines the tumor volume(s) which should receive a specified amount of radiation dose and the critical structures whose dose must not exceed specified tolerance levels. The clinical target volume (CTV) is an extension of the visible tumor in order to account for microscopic spread of cancer cells. 
In order to account for external and internal movements a safety margin around the CTV is applied, which creates the planning target volume (PTV). The organs at risk (OAR) are located asymmetrically around the target creating a minor overlap with the PTV. We follow van Herk, et al.~\cite{vanHerk20001121} who derived the following CTV-PTV-margin-formula 
\begin{equation}\label{eq:vanHerk}
m = 1.96 \Sigma_{tot} + 0.7 \sigma_{tot}
\end{equation}
for a one-dimensional geometry, which satisfies the ICRU criterion \cite{ICRU62} that 99\% of the CTV should receive at least 95\% of prescription dose.
The margin $m$ is derived from the standard deviation of the overall systematic and random errors, $\Sigma_{tot}$ and $\sigma_{tot}$ respectively. 
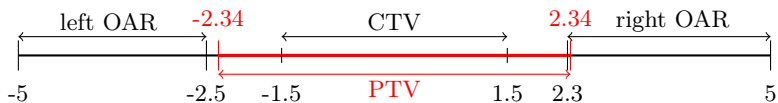
\begin{figure}
\centering
\begin{tikzpicture}
\draw [thick](-5,0) -- (5,0);
\node[draw=none] at (-5,-0.5) {-5};
\node[draw=none] at (5,-0.5) {5};
\node[draw=none] at (-2.5,-0.5) {-2.5};
\node[draw=none,red] at (-2.34,0.5) {-2.34};
\node[draw=none,red] at (2.34,0.5) {2.34};
\node[draw=none] at (2.3,-0.5) {2.3};
\node[draw=none] at (-1.5,-0.5) {-1.5};
\node[draw=none] at (1.5,-0.5) {1.5};
\node[draw=none] at (-3.8,0.45) {\small left OAR};
\node[draw=none] at (0,0.45) {\small CTV};
\node[draw=none] at (3.7,0.45) {\small right OAR};
\node[draw=none][red] at (0,-0.45) {\small PTV};
\draw [very thick,red](-2.34,0) -- (2.34,0);
\draw (-5,-0.2) -- (-5,0.2); 
\draw (-2.5,-0.2) -- (-2.5,0.2); 
\draw [red](-2.34,-0.2) -- (-2.34,0.2); 
\draw (-1.5,-0.1) -- (-1.5,0.1);
\draw (1.5,-0.1) -- (1.5,0.1); 
\draw [red](2.34,-0.2) -- (2.34,0.2); 
\draw (2.3,-0.2) -- (2.3,0.2); 
\draw (5,-0.2) -- (5,0.2); 
\draw [<->] (-5,0.25) -- (-2.5,0.25);
\draw [<->](2.3,0.25) -- (5,0.25);
\draw [<->] (-1.5,0.25) -- (1.5,0.25);
\draw [<->][red] (-2.34,-0.25) -- (2.34,-0.25);
\end{tikzpicture}
\caption{One-dimensional patient phantom. The dimensions are given in cm}\label{fig:Geometry}
\end{figure}

\section{Mathematical Model}
In order to study the concept of adaptive radiation therapy strategies and reduce the computational effort, interfractional variations are simulated by rigid shifts of the whole geometry. By doing so, the movements of each voxel and the dose received at each fraction and in the end of the treatment can be tracked in a straightforward manner. The proposed methods are all designed to take into account  systematic and random geometric errors. The uncertainties may originate for example from patient setup errors, organ motion or inaccuracies in dose delivery systems. It is assumed that random errors differ from fraction to fraction, while the systematic error remains the same over the whole treatment which is divided into $N$ fractions. Therefore, the model assumes the random errors to be independent and identically distributed over the whole course of a treatment. For the purpose of tractability, the geometric errors are accounted for during the optimization process by discretizing their probability distributions into scenarios~\cite{Albin,chen2012including}. The set~$\mathcal{S}$ contains all scenarios of the systematic errors with a the \textit{a priori} normal distribution distribution $p$ with a given mean~$\mu$ and standard deviation~$\Sigma$.  These discrete scenarios are obtained by adapting the probability distribution of the shifts to the voxel-grid according to the mean and standard deviation. The set~$\mathcal{T}$ contains $N$-tuples of the random error scenarios which are assumed to follow the \textit{a priori} normal distribution~$q$ with a given mean~$\mu_q$ and standard deviation~$\sigma$.  Every element $T^N_n$ in an $N$-tuple represents a realized uncertainty scenario for a fraction $n$. In order to reduce dose heterogeneity we introduce the set $\mathcal{U}$ which contains scenarios for  errors following the distribution~$p_{\sigma}$ conditioned on a random standard deviation~$\sigma_u$~\cite{Albin, chan2006robust}. 
\subsection{Nominal Plan}\label{sec:Nominal}
The \textit{nominal plan} corresponds to the procedure of adding safety margins around the CTV. The planned dose is supposed to cover the whole PTV which ensures sufficient CTV coverage during the whole course of the treatment. In the traditional treatment approach the same treatment plan is delivered during all fractions. The optimization problem with the optimization function $f_r$ with $r \in \mathcal{R}$, is formulated according to
\begin{equation}
\begin{aligned} \label{eq:noMCO}
& \underset{x \in \mathcal{X}}{\text{minimize}}
& & \sum_{r \in \mathcal{R}} w_r  f_r(d). \\
\end{aligned}
\end{equation}
The set $\mathcal{R}$ contains the regions of interest, \textit{PTV}, \textit{right OAR}, \textit{left OAR} and \textit{external}, included in the optimization problem, where the importance weights for the corresponding regions are denoted by $w_r$. In this objective function one weighted sum with one set of weights is considered. There is the possibility of extending problem~(\ref{eq:noMCO}) to a multicriteria optimization problem, which is not studied here. The set $\mathcal{X}$ is a collection of all feasible treatment parameters. The dose vector~$d$ is a function of the fluence $x$. The form of the optimization function~$f_r(d)$ is formulated as a quadratic penalty function
\begin{equation} \label{eq:ramp}
f_r(d) = \sum_{v \in \mathcal{V}_r} \Delta_{v,r}(d_v - \hat{d_r})^2,
\end{equation}
where the relative volume of voxel $v$ in volume $r$ is $\Delta_{v,r}$ which satisfies $\sum_{v \in \mathcal{V}_r} \Delta_{v,r} = 1$. The set of voxels for the corresponding structure $r \in \mathcal{R}$ is referred to as $\mathcal{V}_r$. The objective function sets a quadratic penalty for deviations of the voxel dose $d_v$ from the prescribed voxel dose $\hat{d_r}$, aiming at a uniform dose distribution inside the target structure.

\subsection{Robust Plan}\label{sec:combo}
The goal of robust optimization is to create plans which take into account random and systematic errors simultaneously. In contrast to the nominal plan where the PTV is considered to be the high dose target, the focus here lies on optimizing the uniform dose in the CTV in presence of geometric errors. This approach may result in a robust fluence profile similar to applying CTV-PTV-margins. Conditional-value-at-risk (CVaR) optimization~\cite{rockafellar2000} enables to continuously scale between probabilistic and worst case scenario optimization. This is realized by maximizing over the worst systematic error scenarios after a stochastic minimax optimization with respect to the $\alpha$ worst random standard deviations for the expected random errors has been performed. The following problem formulation has been proven to be effective by Fredriksson~\cite{Albin},
\begin{equation}\label{eq:systrand}
\begin{aligned}
& \underset{x \in \mathcal{X}}{\text{minimize}}
&& \underset{\substack{0 \leq \pi \leq \frac{1}{\alpha}p \\ e^T\pi = 1}} {\text{max}} \ \mathbb{E_{\pi}} \bigg[
&& \underset{\substack{0 \leq \rho \leq \frac{1}{\alpha}p_{\sigma} \\ e^T\rho = 1}} {\text{max}} \ \mathbb{E_{\rho}} \big[\mathbb{E_{\textit{q} (\sigma_U)}}[(f(d(x,\mathcal{S},\mathcal{U},\mathcal{T}))]\big] \bigg].
\end{aligned}
\end{equation}
Besides the fluence $x$, the other two variables in the above stated optimization problem~( \ref{eq:systrand}) are the probabilities of the systematic error scenarios $p$ and $p_{\sigma}$ of the unknown random error scenarios, denoted by  $\pi$ and $\rho$ respectively. The closer $\alpha$ is to 1, the less conservative the radiation plan will be. Thus for $\alpha = 1$, the resulting optimization problem~(\ref{eq:systrand}) corresponds to probabilistic optimization~\cite{Albin}. On the other hand, lowering $\alpha$ to $\alpha_s \leq \text{min}_{s \in \mathcal{S}}\  p_s$ with $p_s > 0$ is equivalent to worst case optimization. The problem~(\ref{eq:systrand}) can be written as one minimization problem taking the dual of its inner maximization problem, thus providing continuous derivatives to the optimization solver. 
For every systematic error $ s \in \mathcal{S}$ additional auxiliary
variables $\bar{\lambda}(s)$ and $\bar{\mu}(s)$ are necessary. The
variable $\bar{\lambda}(s)$ is a scalar for the CVaR optimization over
the $\alpha$ worst systematic error scenarios, while $\bar{\mu}(s)$ is
a vector of length~$|\mathcal{U}|$. The optimization problem takes the form
\begin{equation}\label{eq:combdual}
\begin{aligned}
& \underset{\substack{x,\lambda,\mu,\bar{\lambda}, \bar{\mu}}} {\text{minimize}} 
& & \lambda + \frac{1}{\alpha} p^T \mu \\[-2mm]
& \text{subject to} && \lambda + \mu_s \geq \bar{\lambda}(s) + \frac{1}{\alpha}p_\sigma^T \bar{\mu}(s), & \quad s \in \mathcal{S},\\
&&& \bar{\lambda}(s) + \bar{\mu}_u(s) \geq \mathbb{E_{\textit{q} (\sigma_\textit{u})}}[f(d(x,s,u,\mathcal{T}))], & \quad s \in \mathcal{S}, u \in \mathcal{U},\\
& & & \bar{\mu}(s) \geq 0, \quad s \in \mathcal{S}, \quad
\mu \geq 0, \quad
x \in \mathcal{X}.
\end{aligned}
\end{equation}

\section{Adaptive Strategies}
The method described in section \ref{sec:combo} is based on anticipated scenarios conditioned on an \textit{a priori} probability distribution. The optimal solution to~(\ref{eq:combdual}) provides a fluence profile $x \in \mathcal{X}$ which is robust against the included uncertainties. This robust plan is delivered in the first $M-1$ fractions. Rigid whole-body movement can be modeled as isocenter shifts for a one-dimensional patient geometry. These isocenter shifts $\Delta r_1,\ \Delta r_2,...,\Delta r_{M-1}$ are measured during the first $M-1$ fractions. The impact of these treatment errors on the dose distribution is assumed to be analyzed after fraction $M-1$. In case of non-negligible deviations from the planned dose distribution, adaptive replanning based on the actual measured errors has to be considered. The adapted fluence profile, also referred to as plan, is then applied in fraction $M$ and the subsequent fractions until another adaptation may be necessary. Therefore, adapted plans may provide sufficient target coverage while giving little dose to the adjacent healthy tissue for the rest of the treatment.
\subsection{Strategy I: Adaptation of the Probability Distribution}\label{sec:adapt prob}
In the first strategy, the uncertainty sets $\mathcal{S}$ and $\mathcal{T}$ are adapted based on the measured isocenter shifts~$\Delta r_1, \Delta r_2,\dots,\\ \Delta r_{M-1}$. These updated uncertainty sets replace the previous uncertainty sets in the robust optimization problem~(\ref{eq:combdual}). After having modified~(\ref{eq:combdual}), the solver generates an adapted plan. The distribution of the isocenter shifts can be adapted in different ways according to the measurements $\Delta r_0,\Delta r_1, \ \Delta r_2,\dots,\Delta r_{M-1}$ , where $\Delta r_0$ denotes the patient position in the planning CT. 

\begin{description}
\item[{Uncertainty Set Adaptation based on Arithmetic Mean Value}] The first approach is based on the arithmetic mean value. According to Unkelbach and Oelfke \cite{unkelbach2004inclusion} the isocenter position~$\Delta i_{M}$ which is the basis for the updated probability distribution of the systematic error may be determined by the mean value of the measured patient position $\Delta i_M = \frac{1}{M} \sum_{\mu=0}^{M-1} \Delta r_\mu$. The updated standard deviation ${\sigma_i}_M$ is estimated by~${\sigma_i}_M^2 = \frac{1}{M} \sum_{\mu=0}^{M-1} (\Delta r_\mu - \Delta i_M)^2$. Since the uncertainty set for the systematic errors is based on the normal distribution, the mean value $\mu$ and standard deviation $\Sigma$ used in the planning stage are replaced by~$\Delta i_M$ and~${\sigma_i}_M$. The adapted uncertainty scenarios for the systematic errors $s'\in \mathcal{S'}$, the adapted uncertainty set, and their assigned probabilities $p'$ are obtained from discretizing the normal distribution~$\mathcal{N}(\Delta i_M,{\sigma_i}_M^2)$ with respect to the voxel grid. This discretization process of the adapted normal distribution is done in the same way as for the \textit{a priori} distribution~$\mathcal{N}(\mu,\Sigma^2)$. The uncertainty set $\mathcal{S}$ in~(\ref{eq:combdual}) is then replaced by $\mathcal{S'}$. In this work, the adaptive model is extended to include the random errors as well. The adaptation of the random error scenarios is modeled similarly to that of the systematic errors. First, the distances of the measured points$\Delta r_1,\Delta r_2, ...,\Delta r_{M-1}$ relative to the updated mean value~$\Delta i_{M}$ of the systematic error distribution are determined. The updated mean value $\Delta j_M$  and standard deviation~${\sigma_j}_M$ of the random error distribution is obtained from~$\Delta j_{M} = \frac{1}{M} \sum_{\mu=0}^{M-1} (|\Delta r_\mu - \Delta i_M |)$ and~${\sigma_j}_M^2 = \frac{1}{M} \sum_{\mu=0}^{M-1} (\Delta r_\mu - \Delta j_M)^2$.

\item[Uncertainty Set Adaptation based on Exponential Smoothing] The second approach is based on exponential smoothing and is inspired by Chan and Mi{\v{s}}i{\'c} \cite{chan2013adaptive} who used this method to compute the updated mean value~$\Delta i_M$. Let~$\hat{r}_{\mu + 1 }$ denote the forecast of the mean value of the systematic errors for the upcoming fractions~$M, M+1, \dots N$. The forecast is obtained from a convex combination of the so far measured shifts~$\Delta r_{\mu}, \ \text{where}  \ \mu = 0,1,\dots , M-1 $. The weights of the past measurements decay exponentially as the observations get older.
\begin{equation}\label{eq:convex combo}
\hat{r}_{\mu + 1 } = \beta \Delta r_{\mu} + (1-\beta) \hat{r}_{\mu}, \quad
 0\leq \beta \leq 1.
\end{equation}
The advantage of using an update scheme as proposed in (\ref{eq:convex combo}) is the possibility to vary the speed of adaptation by
assigning higher or lower values to the smoothing parameter
$\beta$. The arithmetic mean on the other hand assigns the same weight
to all measurements, which gives a convex combination of the previously computed mean value $\Delta i_{M-1} $ and the most recent measurement $\Delta r_{M-1}$ multiplied by $\frac{1}{M-1}$.
In order to determine the standard deviation $\sigma_{\mu +1}$ for the exponential smoothing method, it is recommended to use the mean absolute value deviation (MAD) and the formula $\sigma \approx (1.25)\text{MAD}$ \cite[chapter 4]{OperationsResearch}. The estimate of MAD for fraction $\mu$ is denoted by $\text{MAD}_\mu$. MAD is the average deviation of each measurement $\Delta r_{\mu}$ from the mean $\hat{r}_{\mu}$. The forecast for MAD in fraction $\mu+1$ is computed by exponential smoothing with the smoothing parameter $\beta$
\begin{equation}
\text{MAD}_{\mu +1}  = \beta | \Delta r_{\mu} - \hat{r}_{\mu}| + (1
-\beta) \text{MAD}_{\mu}, \quad
\hat{\sigma}_{\mu +1} \approx (1.25)\text{MAD}_{\mu +1}.
\label{eq:SDexpSmooth}
\end{equation}
We set $M = \mu + 1$, since we are interested in the forecast for the upcoming fraction $M$. The updated uncertainty set $\mathcal{S}'$ contains scenarios based on the probability distribution $\mathcal{N}( \hat{r}_{M},\hat{\sigma}_{M})$ which is part of the basis for the adapted plan to be delivered in fraction $M$. The next step is to determine the updated set of random errors~$\mathcal{T}'$ using the exponential smoothing method, which is an extension of previous models. Let $\delta_{\mu} = |\Delta r_{\mu} - \hat{r}_{M}|$ denote the deviation of the measurements $\Delta r_0,\Delta r_1,\Delta r_2, ...,\Delta r_{M-1}$  from the recently obtained mean $\hat{r}_{M}$. The computation of the mean random error during the first $M-1$ fractions is carried out according to $\hat{\delta}_{\mu +1} = \beta \delta_{\mu} + (1-\beta) \hat{\delta}_{\mu}  \ \text{with} \ 0\leq \beta \leq 1 $.
The  standard deviation of the random error $\hat{\rho}_M = \hat{\rho}_{\mu +1}$ is computed in the same manner as for the systematic error. The random error scenarios are then based on the discretized normal distribution $\mathcal{N}( \hat{\delta}_{M},\hat{\rho}_{M})$, which are contained in the updated set $\mathcal{T}'$.
\end{description}
The fact that the adapted problem is a new optimization problem is reflected in the notation of the parameters and the objective function. The optimization parameters $ p \ \text{and} \ d $ and the sets $ \mathcal{S} \ \text{and} \ \mathcal{T}$ in (\ref{eq:combdual}) are then replaced by $p', \ d'$, and the sets $ \mathcal{S'} \ \text{and} \ \mathcal{T'}$. 
The solution of this modified robust optimization problem with respect to the adapted uncertainty sets represents the adapted plan to be delivered in fraction $M$ and the remaining $N-M+1$ fractions or until the next adaptation is necessary. 

\subsection{Strategy II: Adaptation of CVaR parameter}\label{sec:adapt CVaR}
The robustness of the solution of the optimization problem (\ref{eq:combdual}) depends on the anticipated uncertainty sets $\mathcal{S}$ and $\mathcal{T}$. The grade of conservativeness is determined by the value that is set for $\alpha$ in the robust problem~(\ref{eq:combdual}). The aim of this strategy is to examine how different values of $\alpha$, i.e.: different levels of conservativeness, are related to the magnitude of systematic and random errors, and how they affect  the cumulative dose distribution. In case of underdosage of the CTV, i.e. large geometric uncertainties, $\alpha$ is going to be decreased in order to improve robustness. In the event of overdosage of a critical structures, $\alpha$ is going to be increased, since a more robust plan may not be necessary to guarantee sufficient CTV coverage in such a scenario. The modified parameter $\alpha'$ replaces $\alpha$ in the objective function and constraint in the original optimization problem~(\ref{eq:combdual}). The solution of this adapted robust problem is an adapted fluence profile to be applied in fraction $M$ and the subsequent $N-M+1$ fractions or until the next adaptation is necessary. 
\subsection{Strategy III: Margin Adaptation}\label{sec:adapt margin}
In the third and last strategy the CTV-PTV-margin is adapted according to the margin formula~(\ref{eq:vanHerk}). The overall distribution of the adapted systematic and random error, $\Sigma'_{tot}$ and $\sigma'_{tot}$ respectively, are based on the obtained measurements $\{ \Delta r_0,\Delta r_1,\Delta r_2, ...,\Delta r_{M-1} \}$. The adapted standard deviations $\Sigma'_{tot}$ and random error $\sigma'_{tot}$ are computed from the arithmetic mean value or by exponential smoothing, as in~(\ref{eq:SDexpSmooth}). The updated margin $m$ is used to define a new PTV which is the target volume in the objective function of the nominal optimization problem~(\ref{eq:noMCO}). This adapted problem is then solved for an updated set of voxels $\mathcal{V}'_r$ and relative volumes $\Delta_{v',r}$ for $r=\text{PTV}$ resulting in an adapted plan for fraction $M$ and the remaining fractions until the next adaptation becomes necessary.

\section{Computational Study}
The computational study is performed in MATLAB version 8.3 using IBM's optimization solver CPLEX in the studio version 12.6.3, where IMRT is applied to the one-dimensional phantom geometry. 

It is assumed that each patient is subjected to one individual systematic error throughout the whole treatment, whereas the random error varies from fraction to fraction. The adaptive strategies and schedules are put to the test by two different patient populations, containing either unpredictably large geometric errors, referred to as \textit{large error population}, or errors with a similar order of magnitude as those accounted for by the initial robust treatment plan. The latter is referred to as \textit{small error population}. In case of the \textit{small error population}, the systematic error is modelled according to the normal distribution around the mean 0 with the standard deviation of 2.5 mm, whereas the random error is distributed around the mean 0 with standard deviation of 6.5 mm. The second population, called the \textit{large error} population, follows a random unknown distribution. The mean of the systematic and random error is assumed to be uniformly distributed $U[-3,3]$. The standard deviation of the systematic error is 3.5 mm and 7.5 mm for the random error.  It should be pointed out that the simulated treatments intentionally follow distributions different from the \textit{a priori} distribution. Simulations show that adaptations are hardly required, if the treatment errors are to similar to the \textit{a prioi} errors. However, one cannot assume that the initial robust plan will be able to handle all types of treatment errors. Therefore, the three proposed strategies are applied to both treatment populations in order to simulate challenging patient cases, which may likely require plan adaptations.
\begin{table}
\center
\begin{tabular}{l|| c | c}
 organ &  relative volume in $\%$ &  relative dose in $\%$\\
\hline
CTV & 99 & 90 \\
left OAR & 20 & 2 \\
right OAR & 30 & 2
\end{tabular}
\caption{Dose-volume points for three structures. The dose levels are defined relative to the prescription dose.}
\label{tab:DVH}
\end{table}
Each of these populations contains 1000 treatments/patients, where the total dose is distributed over 30 fractions for each treatment/patient. At a scheduled evaluation, the so far delivered dose and the dose statistics in the target and critical structures are analyzed. If the measured dose statistics do not satisfy the predefined quality criteria, as listed in Table~\ref{tab:DVH}, the initial plan is adapted. These dose-volume points are recommended by ICRU (International commission on radiation units and measurements) and RTOG (radiation therapy oncology group). This prescription is intended for a prostate cancer case with the prescribed target dose of 70 cGy in presence of the bladder and rectum as critical structures. 

In addition, the impact of adaptation frequency on the adapted plans and cumulative dose distribution is studied by combining the three adaptive strategies with various evaluation schedules. The schedules on trial vary in the total number of evaluations and when the evaluation takes place. They are denoted by \textit{WiEvalj}, where \textit{i} indicates after how many weeks the first evaluation is carried out, whereas the total number of evaluations is given by \textit{j}. 
The search for the right time to evaluate the so far delivered dose and to adapt for the first time during the treatment, is addressed by the schedules \textit{W1Eval4}, \textit{W2Eval4}, \textit{W1Eval3}, \textit{W2Eval3}, \textit{W1Eval1} and \textit{W2Eval1}. The first evaluation is intended after the first and second week of the treatment for the above mentioned schedules. The time interval between each evaluation is one week, with the exception of two weeks between the first and second evaluation for the schedules \textit{W1Eval4} and \textit{W1Eval3}. Moreover, by studying the results of the schedules \textit{W3Eval1}, \textit{W4Eval1} and \textit{W5Eval1} it is expected to estimate after which fraction it is too late to correct for undesired overdosage in the critical structures and underdosage of the CTV. The schedule called \textit{Gold}, refers to the approach of evaluating the to be delivered dose before every fraction followed by a plan adaptation if necessary.
\begin{description}
\item[\textbf{Strategy I}] The robust fluence profile, also referred to as plan, is obtained from solving~(\ref{eq:combdual}) for the parameter $\alpha = 0.1$. This plan has a high level of robustness and is applied during the first fractions, as there is no or limited knowledge of the patient's individual geometric uncertainty distribution. In general, the main objective is to damage the tumor and therefore sufficient target coverage has the highest priority. If the the extent of geometric variations is underestimated, the loss of sufficient target coverage is difficult to compensate for in later fractions~\cite{unkelbach2004inclusion}. If an adaptation is necessary, an updated mean value and standard deviation are computed from the measured errors. The exponential smoothing method is used for a set of three different smoothing parameters $\beta = \left\lbrace 0.1, 0.4, 0.9 \right\rbrace$. Stronger emphasis is put on the most recent information by using a large smoothing parameter, which can be understood as a higher adaptation speed. On the other hand, when multiplying the most recent measurement with a smaller smoothing parameter, the speed of adaptation is lower. Therefore, finding the smoothing parameter best suited to predict the actual individual error distribution is of interest.
\item[\textbf{Strategy II}] The goal here is to determine the appropriate degree of robustness for an adaptive treatment approach. Therefore, the treatments are simulated for initial plans with various degrees of robustness, i.e.: $\alpha$ is either 0.1, 0.4 or 0.9. If the OAR criteria is violated, $\alpha$ will be increased by $0.09$. In case of underdosage of the CTV, the value of the parameter $\alpha$ will decrease by $0.09$. This gradual change in $\alpha$ is carried out in order to develop some proper understanding how different values of $\alpha$ relate to the varying size of geometric uncertainties. 
\item[\textbf{Strategy III}] In this strategy, the nominal plan is applied until the first scheduled evaluation. If one of the predefined quality criteria is not met, an adapted CTV-PTV-margin is used to compute a reoptimized plan according to~(\ref{eq:noMCO}) and~(\ref{eq:ramp}). In contrast to this strategy, the previously described strategies are considered less conservative. 
\end{description}
\section{Results}
\subsection{Non-Adaptive Strategies}
Non-adaptive strategies refer to the traditional approach of using the same plan in all fractions of the treatment. This treatment strategy is intended as a benchmark for further comparison with the adaptive strategies and their potential to improve the overall treatment outcome. 

The nominal and robust plan are applied to the \textit{small error} and \textit{large error population} and evaluated in terms of target coverage and OAR protection, which is illustrated in Table~\ref{tab:CompareSuccessSmallLarge}. The success rate refers to percentage of patients who receive a treatment with a cumulative dose distribution that satisfies all three predefined plan quality criteria. The terms \textit{lowest D99} for the CTV, and \textit{highest D30} and \textit{D20} for both OARs describe the worst possible outcome of the predefined quality criteria. According to the simulations with both populations, higher success rate and much lower dose to the critical structures is achieved by the robust plan. The dose received by the OARs in the worst possible outcomes for both plans is lower in the \textit{large error population} than in the \textit{small error population}. This difference is caused by the wider spread of the geometric errors which smear out the cumulative dose distribution to a greater extent than the errors in the \textit{small error population}. The outstanding performance of the robust-optimized plan applied to \textit{small error distribution} with a success rate of 98\% suggests that adaptation may be only necessary in some exceptional cases, as the \textit{a priori} uncertainty distribution is similar to the simulated treatment uncertainty distribution. A treatment with the nominal plan achieves a lower success rate of 59\% because of too high doses received by the OARs. This difference is not surprising as the robust plan is optimized with respect to the target and the OARs, while it is the nominal plan's main objective to secure target coverage. On the contrary, the simulations with the \textit{large error population} give lower success rates of 2\% and 32\% of the nominal and robust plan, respectively. Thus, combining robust treatment planning with adaptive strategies may increase the probability to satisfy the predefined quality criteria. 

\begin{table}
\center
\begin{tabular}{l | c | c | c | c }
\multirow{2}{*}{Strategy } & success rate in \% &  lowest D99 CTV   & highest D30  right OAR & highest D20 left OAR \\
&  &  (failed patients) in \% & (failed patients) in \% & (failed patients) in \% \\ 
\hline
\multicolumn{5}{l}{ \textit{small error population} }\\
\hline
\hline
Nominal Plan & 59 & 98 (0) & 41 (9) & 82 (14) \\
Robust Plan & 98 & 82 (2) & 4 (0.4) & 16 (0.7)\\
\hline
\multicolumn{5}{l}{ \textit{large error population} }\\
\hline
\hline
Nominal Plan & 2 & 97 (0) & 22 (69) & 35 (82) \\
Robust Plan & 32 & 84 (64) & 8 (16)& 9 (18)\\
\hline
\end{tabular}
\caption{Performance overview of the treatment strategies without any adaptation under trial for the simulated \textit{small} and \textit{large error patient population}. The worst possible outcome in terms of target coverage and organ at risk protection and the percentage of patients failing the specific goal are listed. }
\label{tab:CompareSuccessSmallLarge}
\end{table}

\subsection{Adaptive Strategies}
The number of patients who require an adapted plan at each evaluation is one of the investigated measures. The ratio of these type of patients is referred to as ART candidates and measured in \%.

\begin{figure}
\begin{center}
\subfigure[Dose-Probability histogram comparing the probability that 99\% of the CTV will receive a certain minimum dose or higher.]{\includegraphics[scale=0.40]{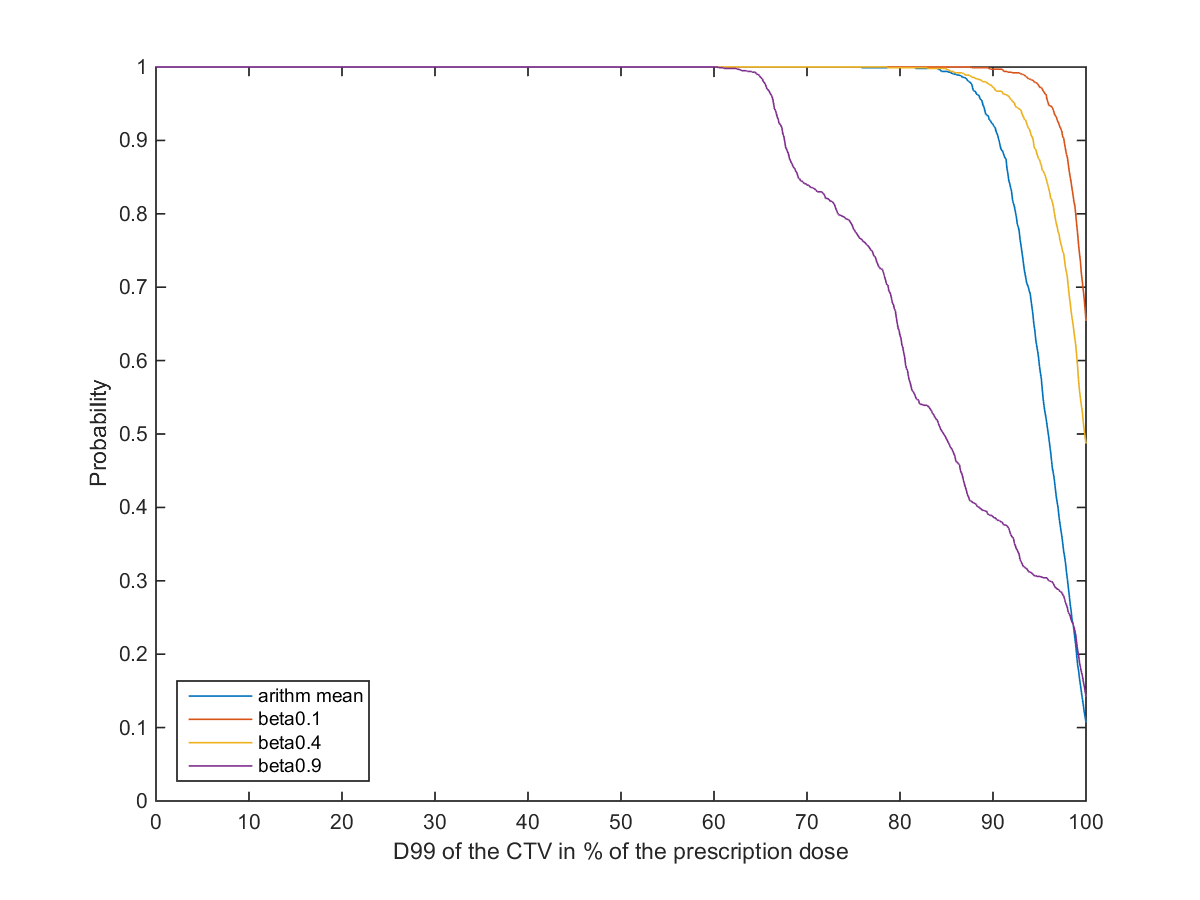}\label{fig:Fig2aCTV}}\quad
\subfigure[Dose-Probability histogram comparing the probability that 30\% of the right OAR will receive a certain maximum dose or lower.]{\includegraphics[scale=0.40]{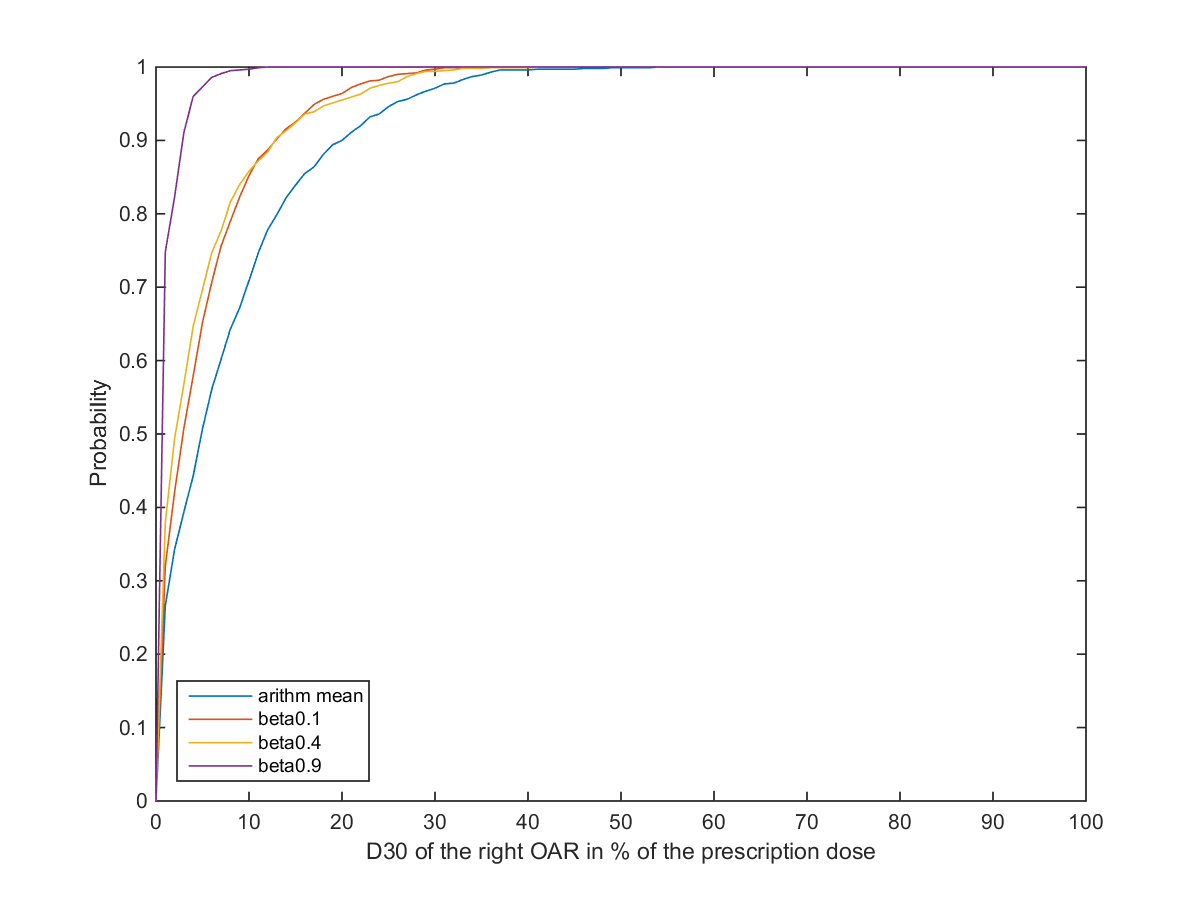} \label{fig:Fig2bOAR}}
\caption{The figure compares four different adaptation speeds for strategy I with the schedule \textit{W1Eval4}, after simulations with the \textit{large error population}.}
\label{fig:doseprobhist}
\end{center}
\end{figure}

\begin{description}
\item[\textbf{Strategy I}] This strategy is applied to both patient populations in order to review the different adaptation schedules in combination with adaptation methods which are based on either arithmetic mean or exponential smoothing. The potential of the adaptive strategy I to improve OAR protection, compared to non-adaptive strategies in presence of unpredictably large variations is similar for the schedules \textit{W1Eval4} and \textit{W2Eval4}. The values of the worst possible outcomes for schedule~\textit{W1Eval4} are illustrated in Table~\ref{tab:PotW1Eval4large}, where the numbers in the parentheses refer to the percentage of patients who fail at the specific goal. Since the right OAR is located closer to the CTV than the left OAR, meeting the predefined quality criterion for the right OAR is challenging. Thus, a higher percentage of patients do not pass this particular quality criterion compared to the number failing the criterion for the left OAR. In general, the quality criteria seem to be satisfied depending on the adaptation method. The highest percentage of succeeding patients for the CTV quality criterion is reached by the exponential smoothing method with $\beta = 0.1$. Furthermore, the heuristic analysis indicates a potential for improved risk organ protection and target coverage, as the overall success rate can be increased to 66\% compared to 32\% after a non-adaptive treatment with the robust plan, as shown in Table~\ref{tab:ResW1andW2Eval4large}. The non-adaptive strategy with the nominal plan gives a success rate of only 2\%.
Concerning the speed of adaptation, the simulations may indicate that the adaptation methods with a higher weight on more recent measurements have difficulties in satisfying the predefined quality criteria at the scheduled evaluation fractions.  The results in Table~\ref{tab:ResW1andW2Eval4large} suggest that strategy I in combination with exponential smoothing does adapt quite accordingly to the simulated uncertainty distribution for a majority of simulated patients.   
By setting a higher weight on earlier measurements, i.e.: a smoothing parameter~$\beta = 0.1, 0.4$, the predictions of future uncertainties may be more accurate than by taking the mean of all measurements~\cite[chapter 4]{OperationsResearch}. The performance of strategy I in presence of errors from the \textit{large error distribution} can be also investigated through dose-probability histograms, as illustrated in Figure~\ref{fig:doseprobhist}. The histograms in Figure~\ref{fig:Fig2aCTV} and~\ref{fig:Fig2bOAR} give information on the probability that the cumulative dose distribution satisfies the predefined quality criteria for the different adaptation methods. The probability, plotted on the y-axis, that 99$\%$ of the CTV will receive at least a certain dose level or above, plotted on the x-axis, is the highest for the exponential smoothing method with~$\beta = 0.1$, as shown in Figure~\ref{fig:Fig2aCTV}. The probability for achieving the same dose levels is slightly lower for the smoothing parameter~$\beta=0.4$, whereas the other two methods reach a much lower probability. On the other hand, setting a greater emphasis on the most recent measurement, as done by setting~$\beta = 0.9$, seems to be most favourable for satisfying the quality criterion for the right OAR, as illustrated in Figure~\ref{fig:Fig2bOAR}. The probability to reach the same dose levels by 30$\%$ of the right OAR is lower for ~$\beta = 0.1 \ \text{and} \ 0.4$. Again, the adaptation method based on the arithmetic mean achieves the lowest probability to fulfil this particular clinical goal. According to the simulations, the adaptation speed has a measurable impact on the probability for sufficient target coverage and OAR protection. A lower adaptation speed, i.e.:~$\beta = 0.1 \ \text{or} \ 0.4$, may lead to a higher probability of target coverage and acceptable OAR protection at the same time. Furthermore, the simulations for the remaining schedules suggest that adaptation is most beneficial, if carried out during the first half of the treatment, as deviations from the intended dose distribution seem to be more difficult to be compensated by later adaptations.
However, applying the same schedules and strategy to the \textit{small error population} suggests that neither the adaptation schedule nor the speed of adaptation significantly affect the final success rate or target coverage. The success rate is comparable to the non-adaptive approach  as the robust plan accounts for treatment errors similar to those in the \textit{small error population}, the success rate of the robust non-adaptive approach is comparable to that of strategy I, independent of the adaptation schedule. Therefore, it can be concluded that not all patients may necessarily have greater benefits from an adaptive approach, if their interfractional motion pattern has been accurately anticipated by the initial robust plan.  

Additional simulations with a subset of 100 cases of the large error population are performed in order to investigate up to which extent the cumulative dose distribution can be improved, if adaptations can be triggered at each fraction. This evaluation schedule under trial is referred to as the \textit{Gold} schedule. The set of simulated treatments is kept small because of long computation times. The simulations with the \textit{large error population} suggest that in presence of unpredictably large geometric uncertainties a higher adaptation frequency may not necessarily lead to an improved success rate, which is illustrated by Table~\ref{tab:GoldLarge}. This could be caused by the high evaluation and adaptation frequency adding noise to the sample of geometric uncertainty measurements. As these measurements build the basis for the reoptimization process, the adapted plans may not be in accordance with the actual uncertainty distribution.

\item[\textbf{Strategy II}] Here, the grade of conservativeness is adjusted in response to unsatisfied predefined dose-volume points. Strategy II is applied to both patient populations in order to understand how robust the initial plan has to be in presence of large or small geometric uncertainties. The simulations indicate that the grade of conservativeness can be lower, i.e.: higher $\alpha$, for a treatment that is designed in an adaptive manner, as the highest success rate is achieved for the value $\alpha = 0.9$ accompanied by a high adaptation rate, for both populations. Throughout the simulated treatments a trend for a decreasing $\alpha$ is observed as the treatment progresses, suggesting that adaptation is triggered by insufficient target coverage. Moreover, it may be concluded that the adjustment of a single parameter may not be sufficient to compensate for previously missed plan-quality criteria, especially in the event of unpredictably large geometric errors from the \textit{large error population}. In case of these larger errors, it is observed that the reoptimized plans tend to underdose the CTV, as the protection of the critical structures might be in conflict with the objective for sufficient target coverage.

\item[\textbf{Strategy III}] The simulations with the \textit{large error population} suggest that the use of CTV-PTV-margins leads to overdosage of the critical structures while ensuring sufficient target coverage. This result is not surprising since target coverage is prioritized in the conventional planning procedure. Due to the broader error distribution the adapted margins increase after the first adaptation and as the treatment progresses. In accordance with previous simulations an adaptive treatment design may not require a conservative initial plan. According to the simulations with the \textit{large error population}, strategy III achieves quite low success rates compared to the previously studied strategies. An increased speed of adaptation, i.e.: large smoothing parameter $\beta = 0.9$, improves the overall success rate from 2\%, in the non-adaptive strategy with the nominal plan, to 21\%. Adaptation of the nominal plan decreases the cumulative dose in the risk organs. Still, the CTV-PTV-margins applied during the treatments may become quite large for patient cases with large geometric variations, as target coverage throughout the treatment is prioritized. Together with the preceding simulations, it may be concluded that a less conservative plan, i.e.: smaller margins, is better suited for an adaptive treatment approach which aims at targeting the CTV and protecting the OARs at the same time. On the other hand, simulations with the \textit{small error population} indicate that strategy~III yields adaptation rates similar to those of the previously tested strategies. 
\end{description}
\begin{table}
\center
\begin{tabular}{l | c | c | c }
\multirow{2}{*}{} &    lowest D99 CTV & highest D30 right OAR & highest D20 left OAR  \\
 & (failed patients) in \% & (failed patients) in \% & (failed patients) in \% \\ 
\hline 
\hline
arithmetic mean &  76 (41) & 54 (65) & 26 (40)  \\
exp. smoothing $\beta = 0.1$  &  88 (3) & 35 (58) & 17 (33)\\
exp. smoothing $\beta = 0.4$  & 79 (13) &  40 (51) & 19 (29)\\
exp. smoothing $\beta = 0.9$  & 60 (70) & 12 (18) & 8 (16)\\
\hline
\end{tabular}
\caption{Comparison of the performance of four different updating procedures for the schedule \textit{W1Eval4} with strategy I subjected to treatment scenarios from the \textit{large error population}. The worst possible outcome in terms of target coverage and organ at risk protection are illustrated and the percentage of patients failing the specific goal are listed in parentheses.}
\label{tab:PotW1Eval4large}
\end{table}
\begin{table}
\center
\begin{tabular}{l c c c c c c }
& \multicolumn{5}{c}{ART candidates in \%} & {success rate in \%}\\
\hline
\hline
Schedule \textit{W1Eval4} & & & \\
fraction number & 5 & 15 & 20 & 25 & $\cap_{\{5,15,20,25\} }$ & 30 \\
\hline
\hline
arithmetic mean & 32 & 60 & 67 & 75 & 32 & 44\\
exp. smoothing $\beta = 0.1$  & 32 & 55 & 58 & 70 & 26 & 66\\
exp. smoothing $\beta = 0.4$  & 32 & 57 & 60 & 66 & 28 & 64 \\
exp. smoothing $\beta = 0.9$  & 32 & 57 & 64 & 71 & 29 & 29 \\
\hline
\hline
Schedule \textit{W2Eval4} & & & \\
fraction number & 10 & 15 & 20 & 25 & $\cap_{\{10,15,20,25\} }$ & 30 \\
\hline
\hline
arithmetic mean & 32 & 58 & 67 & 74 & 22 & 45  \\
exp. smoothing $\beta = 0.1$  & 32 & 54 & 58 & 66 & 26 & 65  \\
exp. smoothing $\beta = 0.4$  & 32 & 57 & 60 & 66 & 27 & 63  \\
exp. smoothing $\beta = 0.9$  & 32 & 58 & 65 & 73 & 30 & 28  \\
\hline
\end{tabular}
\caption{Comparison of the performance of four different updating procedures for the schedules \textit{W1Eval4} and \textit{W2Eval4} with strategy I after simulations with the \textit{large error population}. The percentage of patients whose cumulative dose distribution satisfy the D99, D20 and D30 dose-volume criteria throughout the treatment can be compared for the two schedules.}
\label{tab:ResW1andW2Eval4large}
\end{table}

\begin{table}
\center
\begin{tabular}{c | c | c | c | c }
most common  & success rate   &  lowest D99 CTV  & highest D30 right OAR  & highest D20 left OAR   \\
number of adaptations  & in \% &  (failed patients) in \% & (failed patients) in \% & (failed patients) in \%   \\ 
\hline 
\hline
19 & 49 & 87 (31) & 20 (84) & 7 (38)\\
\hline
\end{tabular}
\caption{Summary of adaptation frequency, -rate, success rate  and worst possible outcomes after simulated treatments of 100 patients from the \textit{large error population} according to the \textit{Gold} schedule. The percentage of patients failing the specific goal is listed in parentheses. The adaptation is performed according to the exponential smoothing method with smoothing parameter~$\beta = 0.1$.}
\label{tab:GoldLarge}
\end{table}
\section{Discussion}

As the scope of this work is to establish a proof-of-concept adaptive framework, non-rigid variations such as tumor deformations are neglected in order to simplify the calculations. Thus, setup and delineation errors are modeled as rigid shifts~\cite{Albin,vanHerk20001121,Lowe2016}. This study is carried out on basis of the results of the heuristic simulations. Due to the simplified treatment- and patient model, the presented numbers for adaptation rates, schedules and strategies should be interpreted and compared in relative terms. The simulations are supposed to give trends how the adaptation of various parameters affect and possibly improve the success rate in the presence of geometric errors. Conclusions are drawn from the treatment simulations with the \textit{large} and \textit{small error population}. Thus, all discussed adaptive and non-adaptive strategies are compared with respect to the same treatment uncertainties and plan quality criteria, which results in a fair comparison of the cumulative dose distributions. In this study, bootstrap resampling is used to asses the statistical accuracy and predictive value of the simulated treatment scenarios. This method is based on the basic idea to randomly draw datasets with replacement from an original dataset. By doing this  $k$ times, a number of $k$ bootstrap datasets are created, where each sample has the same size as the original dataset. Resampling is a helpful tool, whenever studies are subjected to a limited number of patient cases~\cite{Peters1993,DAmico2003}. The largest impact on the dose distribution is caused by the largest geometric errors.  Therefore, the resampling procedure is carried out with the largest errors over the whole treatment of each patient. This procedure is executed over the whole population where 500 datasets are created and the probability density estimate for the standard deviation of the worst case error over the whole large error patient population is computed. As a result, a narrow spread of the standard error is determined, which indicates that the used data set is large enough for quantitative conclusions on the very different performances of the three adaptive strategies. The simulated treatment scenarios and patient geometry turned out to be quite challenging for the robust adaptive strategies to satisfy the predefined quality criteria, as intended from the beginning of the study. Still, it was possible to identify which model parameters and strategy perform well in these challenging treatment scenarios. Therefore, the simulations may represent a guide for further development of robust adaptive treatment strategies. 



\section{Conclusion}
A framework is presented that provides robust optimized and margin adaptation of the treatment plan throughout the fractioned schedule. Measures are introduced which identify patients in need for an adaptive treatment strategy. According to the simulations the combination of evaluation and adaptation schedules with robust reoptimization with respect to predefined plan quality measures may be a useful framework for patient identification. Such treatment strategies open the opportunity to provide individualized plans. Smaller geometric variations are less sensitive to the choice of adaptation strategy as similar geometric variations are already taken into account by the initial robust plan. Still, the dose to the critical structures can be decreased if any of the three strategies is employed. The impact of geometric variations on the dose distribution could be detected and candidates for an adaptive treatment identified.  The lower adaptation rate in the \textit{small error population} is in line with previous studies~\cite{olteanu2014comparative, schwartz2013adaptive} which stated that not all patients may benefit from an adaptive approach in the same extent. The higher adaptation rate for the \textit{large error population}, i.e.: higher number of adaptation candidates, is caused by its larger errors than those included in the \textit{a priori} distribution. However, these simulations show that the cumulative dose distribution can be corrected for a majority of simulated patient treatments. Strategy I seems to achieve the highest success rate and probability for sufficient target coverage and organ at risk protection. Most promising improvements are made when the adaptations are carried out on a weekly basis with the first adaptation done during the first half of the treatment.  Thus, the simulations imply that combining robust treatment planning with adapting the uncertainty scenarios and their assigned probabilities lead to measurable improvements of treatment quality in comparison to conventional non-adaptive strategies. This study is considered to provide a useful base for further development towards more complicated geometries and an optimization model utilizing finite horizon optimal control. 


\begin{acknowledgments}
The authors thank Tatjana Pavlenko for helpful discussion on statistical methods.
\end{acknowledgments}

\bibliography{robust1dliterature}

\begin{thebibliography}{10}

\bibitem{unkelbach2004inclusion}
J.~Unkelbach and U.~Oelfke.
\newblock {Inclusion of organ movements in {IMRT} treatment planning via
  inverse planning based on probability distributions}.
\newblock {\em Physics in Medicine and Biology}, 49(17):4005, 2004.

\bibitem{Albin}
A.~Fredriksson.
\newblock {A characterization of robust radiation therapy treatment planning
  methods-from expected value to worst case optimization}.
\newblock {\em Medical Physics}, 39(8):5169--5181, 2012.

\bibitem{olteanu2014comparative}
L.~Olteanu, D.~Berwouts, I.~Madani, W.~De~Gersem, T.~Vercauteren, F.~Duprez,
  and W.~De~Neve.
\newblock {Comparative dosimetry of three-phase adaptive and non-adaptive
  dose-painting {IMRT} for head-and-neck cancer}.
\newblock {\em Radiotherapy and Oncology}, 111(3):348--353, 2014.

\bibitem{schwartz2013adaptive}
D.~L. Schwartz, A.~S. Garden, S.~J. Shah, G.~Chronowski, S.~Sejpal, D.~I.
  Rosenthal, Y.~I. Chen, Y.~Zhang, L.~Zhang, and P.~F. Wong.
\newblock Adaptive radiotherapy for head and neck cancer - dosimetric results
  from a prospective clinical trial.
\newblock {\em Radiotherapy and Oncology}, 106(1):80--84, 2013.

\bibitem{vanHerk20001121}
M.~van Herk, P.~Remeijer, C.~Rasch, and J.~V. Lebesque.
\newblock The probability of correct target dosage: dose-population histograms
  for deriving treatment margins in radiotherapy.
\newblock {\em International Journal of Radiation Oncology Biology Physics},
  47(4):1121 -- 1135, 2000.

\bibitem{ICRU62}
International~Commission on~Radiation~Units and Measurements.
\newblock {ICRU} report 62: Prescribing, recording, and reporting photon beam
  therapy (supplement to {ICRU} report 50).
\newblock {\em ICRU Publications, Bethesda(MD)}, 1999.

\bibitem{chen2012including}
W.~Chen, J.~Unkelbach, A.~Trofimov, T.~Madden, H.~Kooy, T.~Bortfeld, and
  D.~Craft.
\newblock Including robustness in multi-criteria optimization for
  intensity-modulated proton therapy.
\newblock {\em Physics in Medicine and Biology}, 57(3):591, 2012.

\bibitem{chan2006robust}
T.~C.~Y. Chan, T.~Bortfeld, and J.~N. Tsitsiklis.
\newblock A robust approach to {IMRT} optimization.
\newblock {\em Physics in Medicine and Biology}, 51(10):2567, 2006.

\bibitem{rockafellar2000}
R.~T. Rockafellar and S.~Uryasev.
\newblock Optimization of conditional value-at-risk.
\newblock {\em Journal of Risk}, 2:21--42, 2000.

\bibitem{chan2013adaptive}
T.~C.~Y. Chan and V.~Mi{\v{s}}i{\'c}.
\newblock Adaptive and robust radiation therapy optimization for lung cancer.
\newblock {\em European Journal of Operational Research}, 231(3):745--756,
  2013.

\bibitem{OperationsResearch}
A~.R. Ravindran.
\newblock {\em {Operations Research Applications}}.
\newblock CRC Press, 2008.

\bibitem{Lowe2016}
M.~Lowe, F.~Albertini, A.~Aitkenhead, A.~J. Lomax, and R.~I. MacKay.
\newblock {Incorporating the effect of fractionation in the evaluation of
  proton plan robustness to setup errors.}
\newblock {\em Physics in medicine and biology}, 61(1):413--29, jan 2016.

\bibitem{Peters1993}
L.~J. Peters, H.~Goepfert, K.~K. Ang, R.~M. Byers, M.~H. Maor,
  O.~Guillamondegui, W.~H. Morrison, R.~S. Weber, A.~S. Garden, R.~A.
  Frankenthaler, M.~J. Oswald, and B.~W. Brown.
\newblock {Evaluation of the dose for postoperative radiation therapy of head
  and neck cancer: First report of a prospective randomized trial}.
\newblock {\em International Journal of Radiation Oncology Biology Physics},
  26(1):3--11, apr 1993.

\bibitem{DAmico2003}
A.~V. D'Amico, J.~W. Moul, P.~R. Carroll, L.~Sun, D.~Lubeck, and M.-H. Chen.
\newblock {Surrogate End Point for Prostate Cancer-Specific Mortality After
  Radical Prostatectomy or Radiation Therapy}.
\newblock {\em JNCI Journal of the National Cancer Institute},
  95(18):1376--1383, sep 2003.

\end{thebibliography}
\bibliographystyle{unsrt}

\end{document}